\documentstyle[amstex,amssymb,12pt,a4,thmsa]{article}

\input{tcilatex}
\begin{document}

\title{Approximate McKean-Vlasov Representations for a class of SPDEs\thanks{%
First draft 1 October 2003, current version 22 December 2004.}}
\author{D. Crisan\thanks{%
Department of Mathematics, Imperial College London, 180 Queen's Gate, London
SW7 2BZ, UK.} \ J. Xiong\thanks{%
Department of Mathematics, University of Tennessee, Knoxville, TN 37996-1300
USA. Research supported partially by NSA.}}
\date{22 December 2004}
\maketitle

\begin{abstract}
The solution $\vartheta =\left( \vartheta _{t}\right) _{t\geq 0}$\ of a
class of linear stochastic partial differential equations is approximated
using Clark's robust representation approach (\cite{c}, \cite{cc}). The
ensuing approximations are shown to coincide with the time marginals of
solutions of a certain McKean-Vlasov type equation. We prove existence and
uniqueness of the solution of the McKean-Vlasov equation. The result leads
to a representation of $\vartheta \;$as a limit of empirical distributions
of systems of equally weighted particles. In particular, the solution of the
Zakai equation and that of the Kushner-Stratonovitch equation (the two main
equations of nonlinear filtering) are shown to be approximated the empirical
distribution of systems of particles that have equal weights (unlike those
presented in \cite{kj1} and \cite{kj2}) and do not require additional
correction procedures (such as those introduced in \cite{dan3}, \cite{dan4}, 
\cite{dmm}, etc).\medskip 

\medskip

\end{abstract}

$\left. {}\right. $\newpage

\section{Introduction}

Let $\left( \Omega ,{\cal F},P\right) $ be a probability space on which we
have defined an $m-$dimen\-sional Wiener process $W=\left( W^{i}\right)
_{i=1}^{m}$. Let $M({\Bbb R}^{d})$ be the space of finite Borel measures
defined on the $d-$dimensional Euclidean space ${\Bbb R}^{d}$ and $\left(
\vartheta _{t}\right) _{t\geq 0}$ be an $M({\Bbb R}^{d})$-valued stochastic
process satisfying the following linear stochastic partial differential
equation 
\begin{equation}
\vartheta _{t}(\varphi )=\vartheta _{0}(\varphi )+\int_{0}^{t}\vartheta
_{s}\left( \alpha _{s}\varphi +L_{s}\varphi \right)
ds+\sum_{j=1}^{m}\int_{0}^{t}\vartheta _{s}\left( \beta _{s}^{j}\varphi
\right) dW_{s}^{j},  \label{SPDE}
\end{equation}
where $L=\left( L_{_{s}}\right) _{s\geq 0},$ $\,L_{s}:C_{b}({\Bbb R}%
^{d})\rightarrow C_{b}({\Bbb R}^{d})$ is the second order elliptic
differential operator 
\begin{equation}
L_{_{s}}\varphi =\frac{1}{2}\sum_{i,j=1}^{d}a_{s}^{ij}\partial
_{x_{i}}\partial _{x_{j}}\varphi +\sum_{i=1}^{d}b_{s}^{i}\partial
_{x_{i}}\varphi ,\,\,\,\,\,\,\,\,  \label{operator}
\end{equation}
and $\varphi $ is a function in the domain of $L,$ $\varphi \in
\bigcap_{s\geq 0}{\cal D}\left( L_{s}\right) $. For simplicity, we will
assume that $\vartheta _{0}\;$is a probability measure. If $\alpha \equiv 0$%
, then (\ref{SPDE}) is called the Zakai or the Duncan-Mortensen-Zakai
equation (cf \cite{d2}, \cite{m1}, \cite{za}). The Zakai equation has been
studied extensively over the last 40 years because of its importance in
non-linear filtering theory (see \cite{liptser}, \cite{Pardoux}): In
non-linear filtering, $\bar{\vartheta}_{t},$ the normalized form of $%
\vartheta _{t}$, is the conditional distribution of a (non-homogeneous)
Markov process $\xi _{t}$ with infinitesimal generator $L$ given the
observation process $W$ which satisfies the evolution equation 
\begin{equation}
W_{t}=\int_{0}^{t}\beta _{s}\left( \xi _{t}\right) ds+V_{t},  \label{obser}
\end{equation}
In (\ref{obser}), $V$ is an $m-$dimensional Wiener process independent of $%
\xi $.\ Within the filtering problem $W\;$is\ not a Brownian motion as it is
assumed in the current set-up. However, it becomes a Brownian motion after a
suitable change of measure (Girsanov transformation). One can prove that $%
\bar{\vartheta}_{t}$\ satisfies the following (non-linear) stochastic
partial differential equation
\begin{eqnarray}
\bar{\vartheta}_{t}(\varphi ) &=&\bar{\vartheta}_{0}(\varphi )+\int_{0}^{t}%
\bar{\vartheta}_{s}\left( \alpha _{s}\varphi +L_{s}\varphi \right) -\bar{%
\vartheta}_{s}\left( \alpha _{s}\right) \bar{\vartheta}_{s}\left( \varphi
\right) ds  \nonumber \\
&&+\sum_{j=1}^{m}\int_{0}^{t}\bar{\vartheta}_{s}\left( \beta _{s}^{j}\varphi
\right) -\bar{\vartheta}_{s}\left( \beta _{s}^{j}\right) \bar{\vartheta}%
_{s}\left( \varphi \right) \left( dW_{s}^{j}-\bar{\vartheta}_{s}\left( \beta
_{s}^{j}\right) ds\right) ,  \label{SPDE2}
\end{eqnarray}

If $\alpha \equiv 0$, then equation (\ref{SPDE2}) is called the
Fujisaki-Kallianpur-Kunita or Kushner-Stratonovitch equation (cf. \cite{fkk}%
, \cite{ku}). In general, neither (\ref{SPDE}) nor (\ref{SPDE2}) have
explicit solutions, though one can approximate them by numerical means. As
expected, there is a wide variety of methods to do this (see, for example, 
\cite{dan} and the references therein). Among them, particle methods seem to
be in many cases quite effective. The starting point for a particle
approximation of $\vartheta _{t}$ is the following Feynman-Ka\v{c}
representation 
\begin{equation}
\vartheta _{t}\left( \varphi \right) =E\left[ \varphi \left( X_{t}\right)
A_{t}\left( X\right) |W\right] ,  \label{condrep}
\end{equation}
where $X=\left( X_{t}\right) _{t\geq 0}$ is a Markov process with
infinitesimal generator $L$ independent of $W$ with initial distribution $%
\bar{\vartheta}_{0}$ and $A_{t}\left( X\right) $ is defined as 
\[
A_{t}\left( X\right) =\exp \left( \int_{0}^{t}\alpha _{s}\left( X_{s}\right)
ds+\sum_{j=1}^{m}\left( \int_{0}^{t}\beta _{s}^{j}\left( X_{s}\right)
dW_{s}^{j}-\frac{1}{2}\int_{0}^{t}\beta _{s}^{j}\left( X_{s}\right)
^{2}ds\right) \right) ,\,\text{\thinspace }t\geq 0.
\]
Since in (\ref{condrep}) $X$ is independent of $W$, one can use the Monte
Carlo method to compute $\vartheta _{t}.$ That is, $\vartheta _{t}$ has the
representation 
\begin{equation}
\vartheta _{t}=\lim_{n\rightarrow \infty }\frac{1}{n}\sum_{i=1}^{n}A_{t}%
\left( X^{i}\right) \delta _{X_{t}^{i}},  \label{kx1rep}
\end{equation}
where $X^{i},$ $i>0$ are independent realizations of $X.$ In other words, $%
\left\{ X_{t}^{i},A_{t}\left( X^{i}\right) \right\} _{i=1}^{\infty }$
satisfies the following system of SDEs 
\begin{equation}
\left\{ 
\begin{array}{ll}
dX_{t}^{i} & =b_{t}\left( X_{t}^{i}\right) dt+\sigma _{t}\left(
X_{t}^{i}\right) dB_{t}^{i} \\ 
dA_{t}\left( X^{i}\right)  & =A_{t}\left( X^{i}\right) \alpha _{t}\left(
X_{t}^{i}\right) dt+\sum_{j=1}^{m}A_{t}\left( X_{i}\right) \beta
_{t}^{j}\left( X_{t}^{i}\right) dW_{t}^{j}
\end{array}
\right.   \label{eq}
\end{equation}
where $B^{i},$ $i>0,$ are mutually independent $d-$dimensional Brownian
motions independent of $W$ and $\sigma _{s}=\left( \sigma _{s}\right)
_{i,j=1}^{d}$ is chosen to satisfy $a_{s}=\sigma _{s}\sigma _{s}^{\top }$ ($%
\sigma _{s}^{\top }$ is the transpose of $\sigma _{s}$). As demonstrated in 
\cite{kj1} and \cite{kj2}, representations of the type (\ref{condrep}) and (%
\ref{kx1rep}) hold true for a far wider class of SPDEs than the one
described by (\ref{SPDE}).

However, the convergence in (\ref{kx1rep}) is very slow. That is because the
variance of the weights $A_{t}\left( X^{i}\right) ,$\ $i>0,\;$increases
exponentially fast with time. The effect is that most of weights decrease to 
$0$ with only a few becoming very large. In order to offset this outcome, a
wealth of methods have been proposed. In filtering theory, the generic name
for such a method is that of a{\em \ particle filter} (\cite{dan2}, \cite
{dan3}, \cite{dan4}, \cite{dan5}, \cite{dan6}, \cite{dmm}, etc.). The
standard remedy is to introduce an additional procedure that removes
particles with small weights and adds additional particles in places where
the existing one have large weights. Put differently, one applies at certain
times a branching procedure by which, each particle will be replaced by a
random number of ``offsprings'' with a mean proportional with its
corresponding weight. This branching procedure is a double edge sword:
applying it too often may actually decrease the rate of convergence (cf. 
\cite{dan7}).

In this paper, we seek a different remedy to the slow convergence of the
Monte Carlo method. Heuristically we seek to keep the weights of the
particles equal without introducing an additional procedure but only by
amending the motion of the particles in a way that will take into account
the state of the entire system. This will be achieved in two steps:

First we will show that $\vartheta _{t}$\ and its normalised version $\bar{%
\vartheta}_{t}\;$admit a robust version following the approach first
introduced by Clark (\cite{c}, see also \cite{cc}). By this we mean that $%
\vartheta _{t}$\ and, respectively, $\bar{\vartheta}_{t}\;$depends
continuously on the generating Brownian path $s\rightarrow W_{s}\left(
\omega \right) ,\;\omega \in \Omega $.\ Hence if we consider an
approximating sequence of paths $s\rightarrow W_{s}^{\delta }\left( \omega
\right) \;$that will converge to $s\rightarrow W_{s}\left( \omega \right) \;$%
as $\delta $ tends to $0,\;$then the corresponding measure valued processes $%
\vartheta _{t}^{\delta }$\ and $\bar{\vartheta}_{t}^{\delta }\;$will
converge to $\vartheta _{t}$\ and, respectively, $\bar{\vartheta}_{t}$ with
a certain rate of convergence $r(t,W_{\cdot }\left( \omega \right) ,\delta )$%
.

The second step is to prove that $\bar{\vartheta}_{t}^{\delta }\,\ $has the
asymptotic representation 
\begin{equation}
\bar{\vartheta}_{t}^{\delta }=\lim_{n\rightarrow \infty }\bar{\vartheta}%
_{t}^{\delta ,n},\,\ \ \ \ \,\ \text{where \ }\bar{\vartheta}_{t}^{\delta
,n}=\frac{1}{n}\sum_{i=1}^{n}\delta _{X_{t}^{i,\delta }},  \label{as1}
\end{equation}
and $X^{i},\,$ $i\geq 0\,$are {\em non-linear diffusions} satisfying the
non-linear SDE 
\begin{equation}
dX_{t}^{i,\delta }=\tilde{b}_{t}^{\delta }\left( \tilde{\vartheta}%
_{t}^{\delta },X_{t}^{i,\delta }\right) dt+\sigma _{t}\left( X_{t}^{i,\delta
}\right) dB_{t}^{i}.  \label{new}
\end{equation}
In (\ref{new}), $\tilde{\vartheta}_{t}^{\delta }\,\ $is the conditional
distribution of $X_{t}^{i,\delta }$, $B^{i}\;$are\ mutually independent
Brownian motions and the function $\tilde{b}_{t}^{\delta }\;$will depend
intrinsically on the chosen path $s\rightarrow W_{s}\left( \omega \right) $.
We will show that $\tilde{\vartheta}_{t}^{\delta ,n}=\bar{\vartheta}%
_{t}^{\delta ,n}\;$and that, by suitably choosing $\ $the parameter $\delta
=\delta (W_{\cdot }\left( \omega \right) ,n)$, we will have
\[
\bar{\vartheta}_{t}=\lim_{n\rightarrow \infty }\bar{\vartheta}_{t}^{\delta
,n},\,\ \ \ \ \,\ \text{where \ }\bar{\vartheta}_{t}^{\delta ,n}=\frac{1}{n}%
\sum_{i=1}^{n}\delta _{X_{t}^{i,\delta (W_{\cdot }\left( \omega \right)
,n)}}.
\]
Hence $\bar{\vartheta}_{t}\;$has an asymptotic representation involving
particles with equal weights. Finally in order to obtain the corresponding
asymptotic representation for $\vartheta _{t}$\thinspace we need to
``unnormalize'' $\bar{\vartheta}_{t}^{n}$, i.e., to attach to each particle
a weight $a_{t}^{n}$ (the same one) which converges to $\vartheta _{t}(1)\;$%
as\ $n$ tends to infinity. For example, one could choose
\[
a_{t}^{n}=\exp \left( \int_{0}^{t}\bar{\vartheta}_{t}^{\delta ,n}\left(
\alpha _{s}\right) ds+\sum_{j=1}^{m}\left( \int_{0}^{t}\bar{\vartheta}%
_{t}^{\delta ,n}\left( \beta _{s}^{j}\right) d\left( W_{s}^{\delta }\right)
^{j}-\frac{1}{2}\int_{0}^{t}\bar{\vartheta}_{t}^{\delta ,n}\left( \beta
_{s}^{j}\right) ^{2}ds\right) \right) .
\]
Hence, we obtain the asymptotic representation $\vartheta _{t}\;$which uses
systems of particles with equal weights.
\[
\vartheta _{t}=\lim_{n\rightarrow \infty }\frac{a_{t}^{n}}{n}%
\sum_{i=1}^{n}\delta _{X_{t}^{i,\delta (W_{\cdot }\left( \omega \right)
,n)}}.
\]

We remark that we are not aware of a similar result even for the simpler
case when (\ref{SPDE}) has no stochastic term, i.e., (\ref{SPDE}) is a
second order elliptic PDE with a zero-order term. In this case the first
step described above is not required: One obtains directly an asymptotic
representation of the solution 
\[
\bar{\vartheta}_{t}=\lim_{n\rightarrow \infty }\frac{1}{n}%
\sum_{i=1}^{n}\delta _{X_{t}^{i}} 
\]
where the $X^{i},\,$ $i\geq 0\,$are non-linear diffusions satisfying a
non-linear SDE of the type (\ref{new}).

\section{The robust representation of $\bar{\protect\vartheta}_{t}$}

In this section we introduce the robust representation formula for $%
\vartheta _{t}$ and, respectively, $\bar{\vartheta}_{t}$. By robustness here
we mean that the dependence of the generating Brownian path $t\rightarrow
W_{t}\left( \omega \right) \,$is continuous. The formula is similar to, and
inspired by, the robust version of the integral representation formula of
nonlinear filtering as presented in \cite{c} and \cite{cc}. At a formal
level, the formula is derived by a process of integration-by-parts applied
to the stochastic integrals that appear in the Feynman-Ka\v{c}
representation (\ref{condrep}). The rigorous justification of the formula is
identical with that of the corresponding result in \cite{cc} and for this
reason we omit it here.

For the existence of the robust representation of $\vartheta _{t}$ and,
respectively, of $\bar{\vartheta}_{t}$ we follow Theorems 1 and 2 in \cite
{cc}. For this we need to impose the following conditions:

{\bf RR1}: For all $j=1,...,d$ we assume that\ $\left( s,x\right)
\rightarrow \beta _{s}^{j}\left( x\right) \in C_{1,2}\left( {\Bbb R}%
_{+}\times {\Bbb R}^{d}\right) ,$hence $s\rightarrow \beta _{s}^{j}\left(
X_{s}\right) \,\ $is a semimartingale (we denote by $\beta ^{j}\left(
X\right) ^{m}\;$its martingale part and by $\beta ^{j}\left( X\right) ^{fv}$
its finite variation part) such that, for all $k>0$ and all $t>0$%
\[
E\left[ \exp \left( k\sum_{j=1}^{d}\int_{0}^{t}d\left\langle \beta
^{j}\left( X\right) ^{m}\right\rangle _{s}\right) \right] <\infty ,\,\,\,E%
\left[ \exp \left( k\sum_{j=1}^{d}\int_{0}^{t}\left| d\beta ^{j}\left(
X\right) ^{fv}\right| \right) \right] <\infty . 
\]

{\bf RR2}:\ For all $t>0\;$we have 
\[
E\left[ \exp 2\int_{0}^{t}\alpha _{s}\left( X_{s}\right) ds\right] <\infty . 
\]

Let $y_{\cdot }=\left\{ y_{s},\,s\geq 0\right\} \in C_{{\Bbb R}^{d}}\left(
[0,\infty )\right) $ be a continuous path ${\Bbb R}^{d}\,\ $and define $%
\Theta ^{y_{\cdot }}=\left( \Theta _{t}^{y_{\cdot }}\right) _{t\geq 0}\,$to
be the following $M\left( {\Bbb R}^{d}\right) $-valued stochastic process ($%
\varphi $ is a bounded Borel measurable function) 
\begin{eqnarray*}
\Theta _{t}^{y_{\cdot }}\left( \varphi \right) &=&E\left[ \varphi \left(
X_{t}\right) \exp \left( \int_{0}^{t}\alpha _{s}\left( X_{s}\right)
ds+\sum_{j=1}^{m}\beta _{t}^{j}\left( X_{s}\right) y_{t}^{j}\right. \right.
\\
&&\left. \hspace{1.75in}\left. -\sum_{j=1}^{m}\left( I^{j}\left( y_{\cdot
}\right) +\frac{1}{2}\int_{0}^{t}\beta _{s}^{j}\left( X_{s}\right)
^{2}ds\right) \right) \right] ,
\end{eqnarray*}
where $X=\left( X_{t}\right) _{t\geq 0}$ is a Markov process with
infinitesimal generator $L$ independent of $W$ with initial distribution $%
\bar{\vartheta}_{0}$ and $I^{j}\left( y_{\cdot }\right) \,$is a version of
the stochastic integral $\int_{0}^{t}y_{s}^{j}d\beta _{s}^{j}\left(
X_{s}\right) $. Let $\bar{\Theta}^{y_{\cdot }}=\left( \bar{\Theta}%
_{t}^{y_{\cdot }}\right) _{t\geq 0}\,$be normalized version of the $M\left( 
{\Bbb R}^{d}\right) $-valued stochastic $\bar{\Theta}^{y_{\cdot }}=\left( 
\bar{\Theta}_{t}^{y_{\cdot }}\right) _{t\geq 0}$%
\[
\bar{\Theta}_{t}^{y_{\cdot }}=\frac{\Theta _{t}^{y_{\cdot }}}{\Theta
_{t}^{y_{\cdot }}\left( {\bf 1}\right) }, 
\]
where ${\bf 1}$ is the constant function ${\bf 1}\left( x\right) =1\,$for
all $x\in {\Bbb R}^{d}\,$($\Theta ^{y_{\cdot }}\left( {\bf 1}\right) \,$is
the mass process associated to $\Theta ^{y_{\cdot }}\,$). Following Theorem
1 in \cite{cc}, under the conditions {\bf RR1} and {\bf RR2}, $\Theta
^{y_{\cdot }}\,$and\thinspace $\bar{\Theta}^{y_{\cdot }}\;$depend
continuously on the path $y_{\cdot }$. Moreover, \ for any \ $R>0,\,$%
\thinspace there exist two constants \ $K=K\left( R,t\right) \,\ $and\ $\bar{%
K}=\bar{K}\left( R,t\right) \,$such that\thinspace 
\begin{equation}
\begin{array}{c}
|\Theta _{t}^{y_{\cdot }^{1}}\left( \varphi \right) -\Theta _{t}^{y_{\cdot
}^{2}}\left( \varphi \right) |\leq K\left| \left| \varphi \right| \right|
\left| \left| y_{\cdot }^{1}-y_{\cdot }^{2}\right| \right| _{t} \\ 
\,|\bar{\Theta}_{t}^{y_{\cdot }^{1}}\left( \varphi \right) -\bar{\Theta}%
_{t}^{y_{\cdot }^{2}}\left( \varphi \right) |\leq \bar{K}\left| \left|
\varphi \right| \right| \left| \left| y_{\cdot }^{1}-y_{\cdot }^{2}\right|
\right| _{t}
\end{array}
\label{l1l2}
\end{equation}
for any bounded measurable function $\varphi \;$($\left| \left| \varphi
\right| \right| \triangleq \sup_{x\in {\Bbb R}^{d}}\left| \varphi \left(
x\right) \right| <\infty $) and for any two paths $y_{\cdot }^{1},y_{\cdot
}^{2}\,$such that $\left| \left| y_{\cdot }^{1}\right| \right| _{t},\left|
\left| y_{\cdot }^{2}\right| \right| _{t}\leq R\,\ $(where $\left| \left|
\cdot \right| \right| _{t}\,$is defined as $\left| \left| \alpha \right|
\right| _{t}\triangleq \max_{1\leq i\leq d}\max_{s\in \left[ 0,t\right]
}\left| \alpha _{s}^{i}\right| $).\ Furthermore, if we use the norm $\left|
\left| \cdot \right| \right| _{w}$ on the set of finite signed measures $%
M^{s}\left( {\Bbb R}^{d}\right) \;$%
\[
\left| \left| \mu \right| \right| _{w}=\sup_{\left\{ \varphi \in C_{b}\left( 
{\Bbb R}^{d}\right) |\left| \left| \varphi \right| \right| \leq 1\right\}
}\left| \mu \left( \varphi \right) \right| , 
\]
then, from (\ref{l1l2}), we deduce that 
\begin{equation}
\left| \left| \Theta _{t}^{y_{\cdot }^{1}}-\Theta _{t}^{y_{\cdot
}^{2}}\right| \right| _{w}\leq K\left| \left| y_{\cdot }^{1}-y_{\cdot
}^{2}\right| \right| _{t},\ \ \left| \left| \bar{\Theta}_{t}^{y_{\cdot
}^{1}}-\bar{\Theta}_{t}^{y_{\cdot }^{2}}\right| \right| _{w}\,\leq \bar{K}%
\left| \left| y_{\cdot }^{1}-y_{\cdot }^{2}\right| \right| _{t}.
\label{l1l3}
\end{equation}

Following Theorem 2 in \cite{cc}, $\Theta ^{W_{\cdot }}\,$and\thinspace $%
\tilde{\Theta}^{W_{\cdot }}\,$are the robust versions of $\vartheta $ and,
respectively, $\bar{\vartheta}$. More precisely for all $t\geq 0,$ $\,$%
\begin{equation}
\Theta _{t}^{W_{\cdot }}\,=\vartheta _{t}\,\,\,\text{and\thinspace
\thinspace }\bar{\Theta}_{t}^{W_{\cdot }}=\bar{\vartheta}_{t},\,P\text{%
-almost surely}  \label{l1l4}
\end{equation}
and the null set can be taken to be independent of $t$ \ since all processes
involved are time-continuous.

{\em The robust representation result also enables us to fix the generating
Brownian path }$s\rightarrow W_{s}\left( \omega \right) \;${\em and proceed
to approximate the corresponding }$\bar{\vartheta}_{t}^{\delta }\;${\em and,
implicitly, }$\bar{\vartheta}_{t}\;${\em for this fixed, but arbitrary, path 
}$s\rightarrow W_{s}\left( \omega \right) ${\em . }

We will replace the generating Brownian path $t\rightarrow W_{t}\left(
\omega \right) $ with a smoother one. For example, we choose a partition of
the time interval $[0,\infty )\,\ $of the form 
\[
0=t_{0}<t_{1}<....<t_{i}<...
\]
where $t_{i}=i\delta ,\;i\geq 0\;$and define the following piecewise linear
path $t\rightarrow W_{t}^{\delta }\left( \omega \right) $ 
\[
W_{t}^{\delta }\left( \omega \right) =W_{t_{i}}\left( \omega \right) +\frac{%
W_{t_{i+1}}\left( \omega \right) -W_{t_{i}}\left( \omega \right) }{\delta }%
\left( t-t_{i}\right) ,\,\,\ \text{for}\,\,t\in \lbrack t_{i},t_{i+1}).
\]
Let $\vartheta _{t}^{\delta }\triangleq \Theta _{t}^{W_{\cdot }^{\delta }}\;$%
and $\bar{\vartheta}_{t}^{\delta }\triangleq \bar{\Theta}_{t}^{W_{\cdot
}^{\delta }}$. Then, from (\ref{l1l3}), it follows that, 
\[
\lim_{\delta \rightarrow 0}\vartheta _{t}^{\delta }=\vartheta _{t}\;\;\text{%
and\thinspace \thinspace }\lim_{\delta \rightarrow 0}\bar{\vartheta}%
_{t}^{\delta }=\bar{\vartheta}_{t},
\]
where the above convergence is in the weak topology on $M\left( {\Bbb R}%
^{d}\right) $. Moreover, we have the following rates of convergence: there
exist two constants $K=K\left( \varpi _{t}\left( \omega \right) ,t\right) $
and $\bar{K}=\bar{K}\left( \varpi _{t}\left( \omega \right) ,t\right) $
where $\varpi _{t}\left( \omega \right) =\sup_{s\in \lbrack 0,t+\delta
]}\max_{1\leq i\leq d}\left| W_{s}^{i}\left( \omega \right) \right| $ $\,$%
such that\thinspace\ 
\[
\left| \left| \vartheta _{t}^{\delta }-\vartheta _{t}\right| \right|
_{w}\leq K\varpi _{t}^{\delta }\left( \omega \right) ,\,\,\ \ \ \left|
\left| \bar{\vartheta}_{t}^{\delta }-\bar{\vartheta}_{t}\right| \right|
_{w}\leq \bar{K}\varpi _{t}^{\delta }\left( \omega \right) 
\]
and $\varpi _{t}^{\delta }\left( \omega \right) =\max_{j=0,\left[ t\delta %
\right] }\sup_{t\in \left[ t_{i},t_{i+1}\right] }\max_{1\leq i\leq d}\left|
W_{t}^{i}\left( \omega \right) -W_{t_{j}}^{i}\left( \omega \right) \right| $%
. Moreover since $s\rightarrow W_{s}\left( \omega \right) $ is continuous,
we can choose $\delta =\delta \left( \omega ,n\right) ,\;$so that $\varpi
_{t}^{\delta }\left( \omega \right) \leq \frac{1}{n}.$ Hence 
\begin{equation}
\begin{array}{c}
\left| \left| \vartheta _{t}^{\delta \left( \omega ,n\right) }-\vartheta
_{t}\right| \right| _{{\cal M}}\leq \left| \left| \vartheta _{t}^{\delta
\left( \omega ,n\right) }-\vartheta _{t}\right| \right| _{w}\leq \frac{K}{%
\sqrt{n}}. \\ 
\left| \left| \bar{\vartheta}_{t}^{\delta \left( \omega ,n\right) }-\bar{%
\vartheta}_{t}\right| \right| _{{\cal M}}\leq \left| \left| \bar{\vartheta}%
_{t}^{\delta \left( \omega ,n\right) }-\bar{\vartheta}_{t}\right| \right|
_{w}\leq \frac{\bar{K}}{\sqrt{n}}.
\end{array}
\label{rate2}
\end{equation}
In (\ref{rate2}), the norm $\left| \left| \cdot \right| \right| _{{\cal M}.}$
\footnote{%
Both norms $\left| \left| \cdot \right| \right| _{w}\;$and .$\left| \left|
\cdot \right| \right| _{{\cal M}.}\;$induce the weak topology on\ $M\left( 
{\Bbb R}^{d}\right) $.} is defined as 
\begin{equation}
\left| \left| \mu \right| \right| _{{\cal M}}\triangleq \sum_{k=0}^{\infty }{%
\frac{|\mu \left( \varphi _{k}\right) |}{2^{k}\Vert \varphi _{k}\Vert }%
,\,\,\,\mu \in }M^{s}\left( {\Bbb R}^{d}\right) ,  \label{mnorm}
\end{equation}
where $\varphi _{0}=1$\ and $(\varphi _{k})_{k>0}$ are the elements of a
dense set ${\cal M}\in C_{b}\left( {\Bbb R}^{d}\right) $. This norm will be
used to prove convergence in the next step of the construction.

It is easy to see that $\vartheta _{t}^{\delta }\;$satisfies the following
linear PDE (written in weak form) 
\begin{equation}
\vartheta _{t}^{\delta }(\varphi )=\vartheta _{0}(\varphi
)+\int_{0}^{t}\vartheta _{s}^{\delta }\left( \alpha _{s}^{\delta }\varphi
+L_{s}\varphi \right) ds  \label{PDE1}
\end{equation}
where 
\begin{equation}
\alpha _{s}^{\delta }=\alpha _{s}-\frac{1}{2}\sum_{j=1}^{m}\left( \left(
\beta _{s}^{j}\right) ^{2}+\beta _{s}^{j}\frac{W_{t_{i+1}}^{i}\left( \omega
\right) -W_{t_{i}}^{i}\left( \omega \right) }{\delta }\right) ,\,\ \ \,s\in
\lbrack t_{i},t_{i+1})  \label{alpha}
\end{equation}
whilst $\bar{\vartheta}_{t}^{\delta }\,$satisfies the following nonlinear
PDE (again, written in weak form) 
\begin{equation}
\bar{\vartheta}_{t}^{\delta }(\varphi )=\bar{\vartheta}_{0}(\varphi
)+\int_{0}^{t}\bar{\vartheta}_{s}^{\delta }\left( \alpha _{s}^{\delta
}\varphi +L_{s}\varphi \right) -\bar{\vartheta}_{s}^{\delta }\left( \alpha
_{s}^{\delta }\right) \bar{\vartheta}_{s}^{\delta }\left( \varphi \right) ds
\label{PDE2}
\end{equation}
Note that, from (\ref{PDE1}), $\vartheta _{t}^{\delta }(1)=\exp \left(
\int_{0}^{t}\bar{\vartheta}_{s}^{\delta }\left( \alpha _{s}^{\delta }\right)
ds\right) $ (since $\vartheta _{0}^{\delta }(1)=1$).

We introduce now additional conditions on the coefficients of the SPDE (\ref
{SPDE})\ to insure it has a unique solution and that the solution has a
density with respect to the Lebesgue measure and its density is strictly
positive everywhere.

{\bf EU:} We will assume that the functions $\alpha _{s}$, $\beta _{s}^{j},$ 
$j=1,...,m,\,$are bounded (with a bound independent of $s$). $%
a_{s}^{ij},i,j=1,...,d,\,\,b_{s}^{i},\,\,i=1,...,d$ are bounded (with a
bound independent of $s$) and Lipschitz (with a Lipschitz constant
independent of $s$).\ We also assume that the operator $L$ is uniformly
elliptic. That is there is a constant $b$ such that for all times $s\geq 0\;$%
and vectors $\xi =\left( \xi _{i}\right) _{i=1}^{d}\in {\Bbb R}^{d}$, we
have 
\[
\left\langle \xi ,a_{s}\xi \right\rangle \geq b\left\langle \xi ,\xi
\right\rangle , 
\]
where $\left\langle \cdot ,\cdot \right\rangle \;$is the standard inner
product on ${\Bbb R}^{d}$. Further, we will assume that $\vartheta _{0}\;$%
has finite second moment.

{\bf PD:}\ We will assume that $\vartheta _{0}\;$is absolutely continuous
with respect to the Lebesgue measure and its density is strictly positive
everywhere.

Assuming condition {\bf EU} ensures that the SPDE (\ref{SPDE}) (see for
example Jie\&Kurtz) and the PDEs (\ref{PDE1}) and (\ref{PDE2}) have unique
solutions. Further, the system of equations (\ref{eq}) and, in particular,
the SDE satisfied by the processes $X_{i}\;$(which appear in the asymptotic
representation (\ref{kx1rep}) of $\vartheta _{t}$ ) has a unique solution.
In particular this ensures the existence of a Markov process $X=\left(
X_{t}\right) _{t\geq 0}$ with infinitesimal generator $L$ and initial
distribution $\bar{\vartheta}_{0}.$ Assuming {\bf EU+PD}$\;$guarantees that,
for any $t\geq 0$ the distribution of $X_{t}\;$is absolutely continuous with
respect to the Lebesgue measure and its density is strictly positive
everywhere.

Let us note that, from the Feynman-Ka\v{c} representation (\ref{condrep}), $%
\vartheta _{t}$ (and hence its normalized version $\bar{\vartheta}_{t}$) has
the same support as the distribution of $X_{t}$. Similarly, $\vartheta
_{t}^{\delta }\;$has the following representation 
\begin{equation}
\vartheta _{t}^{\delta }\left( \varphi \right) =E\left[ \varphi \left(
X_{t}\right) \exp \left( \int_{0}^{t}\alpha _{s}^{\delta }\left(
X_{s}\right) ds\right) \right] .  \label{fkp}
\end{equation}
So, both $\vartheta _{t}^{\delta }$ and $\bar{\vartheta}_{t}^{\delta }$ have
the same support as the distribution of $X_{t},$ too. It follows from the
above that $\vartheta _{t}$, $\bar{\vartheta}_{t}$, $\vartheta _{t}^{\delta
} $ and $\bar{\vartheta}_{t}^{\delta }\;$are all absolutely continuous with
respect to the Lebesgue measure and their density is strictly positive
everywhere. In the following we will denote by $x\rightarrow \bar{\vartheta}%
_{t}\left( x\right) \;$and $x\rightarrow \bar{\vartheta}_{t}^{\delta }\left(
x\right) $ the density of $\bar{\vartheta}_{t}$ and respectively $\bar{%
\vartheta}_{t}^{\delta }\;$with respect to the Lebesgue measure.

\section{The non-linear process}

We introduce next the non-linear SDEs satisfied by the non-linear diffusions 
$X^{i,\delta },\,$ $i\geq 0\;$appearing in the asymptotic representation (%
\ref{as1}) of $\bar{\vartheta}_{t}^{\delta }$. For this we need to define
first the drift coefficient of the non-linear SDE (\ref{new}). In the
following, we fix the generating Brownian path $t\rightarrow W_{t}\left(
\omega \right) .$\ The coefficient $\alpha _{s}^{\delta }\;$obviously
depends on the path, however we will not make the dependence explicit. For
this fixed path, the non-linear diffusions $X^{i,\delta }$  will be defined
on a new probability space $\left( \tilde{\Omega},{\cal \tilde{F}},\tilde{P}%
\right) $.

For arbitrary $f\in {\cal B}\left( {\Bbb R}^{d}\right) $ and $\mu \in {\cal P%
}\left( {\Bbb R}^{d}\right) $, define $f^{\mu }\in {\cal B}\left( {\Bbb R}%
^{d}\right) $ to be the function $f^{\mu }\triangleq f-\mu \left( f\right) $
and let $\Lambda _{\mu }f=\left( \Lambda _{\mu }^{j}f\right) _{j=1}^{d}$ be
the vector function 
\[
\Lambda _{\mu }f\left( x\right) =\frac{1}{\omega _{d}}\int_{{\Bbb R}^{d}}%
\frac{\left( y-x\right) f^{\mu }\left( y\right) }{\left| \left| x-y\right|
\right| ^{d}}\mu \left( dy\right) ,
\]
where $\omega _{d}\,\ $is the surface area of the \ $d-$dimensional sphere \ 
$S_{d-1}$. If $d=1,$ then $\Lambda _{\mu }f$ has several equivalent
representions representation 
\begin{eqnarray}
\Lambda _{\mu }f\left( x\right)  &=&\mu \left( f1_{[x,\infty )}\right) -\mu
\left( f\right) \mu \left( 1_{[x,\infty )}\right) =\mu \left( f^{\mu
}1_{[x,\infty )}\right)   \nonumber \\
&=&\mu \left( 1_{\left( -\infty ,x\right) }\right) \mu \left( f\right) -\mu
\left( f1_{\left( -\infty ,x\right) }\right) =-\mu \left( f^{\mu }1_{\left(
-\infty ,x\right) }\right)   \label{onerep}
\end{eqnarray}
When $d=1$ the function $\Lambda _{\mu }f$ is well defined for arbitrary \ $%
f\in {\cal B}\left( {\Bbb R}^{d}\right) $ and $\mu \in {\cal P}\left( {\Bbb R%
}^{d}\right) $. This is not the case when\ $d>2$. Let $\mu \in {\cal P}%
\left( {\Bbb R}^{d}\right) ,d\geq 2$ be a probability measure absolutely
continuous with respect to the Lebesgue measure. We say that $\mu $ is a $p-$%
{\em good} probability measure if its density $x\rightarrow d_{\mu }\left(
x\right) $ with respect to the Lebesgue measure is locally bounded and \ $p$%
-integrable, where \ \ $1<p<d$. The following proposition gives a necessary
condition which insures that $\Lambda _{\mu }f$ is well defined:

\begin{proposition}
\label{pwell}If $\mu \in {\cal P}\left( {\Bbb R}^{d}\right) ,d\geq 2$ is a $%
p $-good probability measure, then $\Lambda _{\mu }f$ is well defined.
\end{proposition}

\proof%
%
We need to prove that the function 
\[
y\longrightarrow \frac{y-x}{\left| \left| x-y\right| \right| ^{d}}\left(
f\left( y\right) -\mu \left( f\right) \right) 
\]
is integrable with respect to $\mu $ for all $x>0$. Let $\left| \left|
f\right| \right| $ be the \ $L^{\infty }$-norm \ of \ $f$ and $\left| \left|
\left. \mu \right| _{B\left( x,1\right) }\right| \right| $ be the $L^{\infty
}$-norm \ of the function $x\rightarrow d_{\mu }\left( x\right) \,\,$\
restricted to $B\left( x,1\right) ,$ the ball of radius \ $1$ and centre \ $x
$. Then 
\begin{eqnarray}
\frac{1}{\omega _{d}}\int_{{\Bbb R}^{d}}\left| \frac{y_{i}-x_{i}}{\left|
\left| x-y\right| \right| ^{d}}\left( f\left( y\right) -\mu \left( f\right)
\right) \right| \mu \left( dy\right)  &\leq &\frac{2\left| \left| f\right|
\right| \left| \left| \left. \mu \right| _{B\left( x,1\right) }\right|
\right| }{\omega _{d}}\int_{B\left( x,1\right) }\frac{1}{\left| \left|
x-y\right| \right| ^{d-1}}dy  \nonumber \\
&&+\int_{B\left( x,1\right) ^{c}}\frac{1}{\left| \left| x-y\right| \right|
^{d-1}}\mu \left( dy\right) ,  \label{welldef}
\end{eqnarray}
where $B\left( x,1\right) ^{c}\triangleq {\Bbb R}^{d}\backslash B\left(
x,1\right) .\;$By using polar coordinates it is easy to check that the first
integral is equal to $\omega _{d}$ and hence it is finite. Using H\"{o}%
lder's inequality and polar coordinates, we get that 
\begin{eqnarray*}
\int_{B\left( x,1\right) ^{c}}\frac{1}{\left| \left| x-y\right| \right|
^{d-1}}\mu \left( dy\right)  &\leq &\left( \int_{B\left( x,1\right) ^{c}}%
\frac{1}{\left| \left| x-y\right| \right| ^{q\left( d-1\right) }}dy\right) ^{%
\frac{1}{q}}\left( \int_{B\left( x,1\right) ^{c}}\left( d_{\mu }\left(
y\right) \right) ^{p}dy\right) ^{\frac{1}{p}} \\
&\leq &\left( \omega _{d}\int_{1}^{\infty }r^{-\left( q-1\right) \left(
d-1\right) }dr\right) ^{\frac{1}{q}}\left| \left| \mu \right| \right| _{p}
\end{eqnarray*}
and since $\left( q-1\right) \left( d-1\right) =\frac{d-1}{p-1}>1$ we get
that $\int_{1}^{\infty }r^{-\left( q-1\right) \left( d-1\right) }dr$ is
finite. Hence the second integral on the right hand side of (\ref{welldef})
is finite and so $\Lambda _{\mu }f$ \ is well defined and is a bounded
function.%
\endproof%
%
\bigskip 

\noindent We now define the first order differential operators 
\[
L_{_{\mu }}^{f}\varphi \left( x\right) \triangleq \sum_{j=1}^{d}\Lambda
_{\mu }^{j}f\left( x\right) \partial _{x_{j}}\varphi \left( x\right)
,\,\,\,\ \bar{L}_{_{\mu }}^{f}\varphi \left( x\right) \triangleq \frac{%
L_{_{\mu }}^{f}\varphi \left( x\right) }{d_{\mu }\left( x\right) }. 
\]
In the following, we will \ denote \ by $l_{d}$ to be the Lebesgue measure
on \ ${\Bbb R}^{d}.$

\begin{proposition}
\label{parts}$\left. {}\right. $\newline
i. If \ $\varphi \in {\cal C}_{b}^{1}\left( {\Bbb R}\right) $, then 
\[
\mu \left( \varphi f^{\mu }\right) =l_{1}\left( L_{_{\mu }}^{f}\varphi
\right) =\mu \left( \bar{L}_{_{\mu }}^{f}\varphi \right) . 
\]
ii. If \ $d\geq 2$, $\varphi \in {\cal C}_{k}^{1}\left( {\Bbb R}^{d}\right) $
and $\mu \in {\cal P}\left( {\Bbb R}^{d}\right) ,d\geq 2$ is a $p$-good
probability measure ($1<p<d$), then 
\[
\mu \left( \varphi f^{\mu }\right) =l_{d}\left( L_{_{\mu }}^{f}\varphi
\right) =\mu \left( \bar{L}_{_{\mu }}^{f}\varphi \right) . 
\]
\end{proposition}

\proof%
%
$\left. {}\right. $\newline
i. Using the representation (\ref{onerep}) we have that $\lim_{\left|
x\right| \rightarrow \infty }\Lambda _{\mu }^{j}f\left( x\right) =0,$ hence
by integration by parts 
\[
\int_{{\Bbb R}}\varphi \left( x\right) f^{\mu }\left( x\right) \mu \left(
dx\right) =\int_{{\Bbb R}}\varphi ^{\prime }\Lambda _{\mu }f\left( x\right)
dx=\int_{{\Bbb R}}L_{_{\mu }}^{f}\varphi \left( x\right) dx=\int_{{\Bbb R}}%
\bar{L}_{_{\mu }}^{f}\varphi \left( x\right) \mu \left( dx\right) \text{.} 
\]
ii. For $\varphi \in {\cal C}_{1}^{k}\left( {\Bbb R}^{d}\right) $ the
following identity holds true (see, for example, \cite{lsc} page \ 12) 
\[
\varphi \left( y\right) =\frac{1}{\omega _{d}}\int_{{\Bbb R}^{d}}\frac{%
\left\langle y-x,\nabla \varphi \left( x\right) \right\rangle }{\left|
\left| x-y\right| \right| ^{d}}dx. 
\]
We will show that the function 
\[
\left( x,y\right) \in {\Bbb R}^{d}\times {\Bbb R}^{d}\longrightarrow \frac{%
\left\langle y-x,\nabla \varphi \left( x\right) \right\rangle }{\left|
\left| x-y\right| \right| ^{d}}f^{\mu }\left( y\right) d_{\mu }\left(
y\right) 
\]
is integrable with respect to $l_{2d}$. Hence by Fubini's theorem 
\begin{eqnarray*}
\mu \left( \varphi f^{\mu }\right) &=&\int_{{\Bbb R}^{d}}\varphi \left(
y\right) f^{\mu }\left( y\right) d_{\mu }\left( y\right) dy \\
&=&\int_{{\Bbb R}^{d}}\frac{1}{\omega _{d}}\int_{{\Bbb R}^{d}}\frac{%
\left\langle y-x,\nabla \varphi \left( x\right) \right\rangle }{\left|
\left| x-y\right| \right| ^{d}}dxf^{\mu }\left( y\right) d_{\mu }\left(
y\right) dy \\
&=&\int_{{\Bbb R}^{d}}\left\langle \frac{1}{\omega _{d}}\int_{{\Bbb R}^{d}}%
\frac{\left( y-x\right) f^{\mu }\left( y\right) }{\left| \left| x-y\right|
\right| ^{d}}d_{\mu }\left( y\right) dy,\nabla \varphi \left( x\right)
\right\rangle \\
&=&\int_{{\Bbb R}^{d}}\left\langle \Lambda _{\mu }f\left( x\right) ,\nabla
\varphi \left( x\right) \right\rangle =\int_{{\Bbb R}^{d}}L_{_{\mu
}}^{f}\varphi \left( x\right) dx.
\end{eqnarray*}
First let us observe that 
\[
\left| \frac{\left\langle y-x,\nabla \varphi \left( x\right) \right\rangle }{%
\left| \left| x-y\right| \right| ^{d}}f^{\mu }\left( y\right) d_{\mu }\left(
y\right) \right| \leq 2\left| \left| f\right| \right| \frac{\left| \left|
\nabla \varphi \left( x\right) \right| \right| }{\left| \left| x-y\right|
\right| ^{d-1}}d_{\mu }\left( y\right) . 
\]
Let now $I_{1}$ be the following Riesz type operation 
\[
I_{1}\psi \left( y\right) =\int_{{\Bbb R}^{d}}\frac{\psi \left( x\right) }{%
\left| \left| x-y\right| \right| ^{d-1}}dx 
\]
for $\psi \in {\cal C}^{k}\left( {\Bbb R}^{d}\right) $. Theorem 1.2.1 from 
\cite{lsc}, page 12, states that $I_{1}\psi $ is \ $q$-integrable for any $q$
such that 
\[
\frac{1}{q}=\frac{1}{p^{\prime }}-\frac{1}{d}\text{ where \ }1<p^{\prime }<d%
\text{.} 
\]
So $I_{1}\psi $ is $q$-integrable for any $q$ such that $0<\frac{1}{q}<1-%
\frac{1}{d}.$ Hence $I_{1}\psi $ is $q$-integrable \ for $q$ such that $%
\frac{1}{q}=1-\frac{1}{p}$. So, by Fubini's theorem (for non-negative
functions) and H\"{o}lder's inequality, 
\begin{eqnarray*}
\int_{{\Bbb R}^{d}\times {\Bbb R}^{d}}\frac{\left| \left| \nabla \varphi
\left( x\right) \right| \right| }{\left| \left| x-y\right| \right| ^{d-1}}%
d_{\mu }\left( y\right) dxdy &=&\int_{{\Bbb R}^{d}}\left( \int_{{\Bbb R}^{d}}%
\frac{\left| \left| \nabla \varphi \left( x\right) \right| \right| }{\left|
\left| x-y\right| \right| ^{d-1}}dx\right) d_{\mu }\left( y\right) dy \\
&=&\int_{{\Bbb R}^{d}}I_{1}\psi \left( y\right) d_{\mu }\left( y\right) dy \\
&\leq &\left( \int_{{\Bbb R}^{d}}I_{1}\psi \left( y\right) ^{q}dy\right) ^{%
\frac{1}{q}}\left( \int_{{\Bbb R}^{d}}d_{\mu }\left( y\right) ^{p}dy\right)
^{\frac{1}{p}}<\infty ,
\end{eqnarray*}
where $\psi \left( x\right) \equiv \left| \left| \nabla \varphi \left(
x\right) \right| \right| $. Hence our claim.%
\endproof%
%
\bigskip

We are now ready to define the coefficients $\tilde{b}_{s}\,$of the equation
(\ref{new}). Let $\mu $ be a probability measure that satisfies the
following conditions:

1. $\mu $ is absolutely continuous with respect to the Lebesgue measure and
its density $x\rightarrow d_{\mu }\left( x\right) $ with respect to the
Lebesgue measure is a strictly positive function in $C_{b}^{1}({\Bbb R}^{d})$%
.

2. If $d\geq 2$, then there exists $p$ such that $\mu $ is a $p-${\em good}
probability measure.

Then the following coefficients are well defined 
\begin{equation}
\left( s,\mu ,x\right) \longrightarrow \tilde{b}_{s}^{\delta }\left( \mu
,x\right) \triangleq b\left( x\right) +\frac{\Lambda _{\mu }\alpha
_{s}^{\delta }\left( x\right) }{d_{\mu }\left( x\right) },j=1,...,m.
\label{c1}
\end{equation}

Let $X=\left( X_{t}\right) _{t\geq 0}$ be a continuous $d-$dimensional
stochastic process. In the following we will denote by $\tilde{\vartheta}%
_{t} $ the distribution of $X_{t}$. We say that \ $X$ is a {\em good}
process if the following three conditions are satisfied:

\begin{itemize}
\item  $\tilde{\vartheta}_{t}$ is absolutely continuous with respect to the
Lebesgue measure and its density $x\rightarrow \tilde{\vartheta}_{t}\left(
x\right) $ is positive for all \ $x\in {\Bbb R}^{d}$.

\item  The function \ $\left( t,x\right) \in \left[ 0,\infty \right] \times $
${\Bbb R}^{d}\longrightarrow \tilde{\vartheta}_{t}\left( x\right) \,\,$is
continuous.

\item  If $d\geq 2$, then there exists $p$ such that $\tilde{\vartheta}_{t}$
is a $p-${\em good} probability measure for all \ $t\geq 0$.
\end{itemize}

We have the following proposition

\begin{proposition}
\label{nonproc}Let $\left( \tilde{\Omega},{\cal \tilde{F}},\tilde{P}\right)
\;$be a probability space on which there exists a good process $X=\left(
X_{t}^{\delta }\right) _{t\geq 0}$ which satisfies equation (\ref{new}).
That is 
\[
X_{t}^{\delta }=X_{0}^{\delta }+\int_{0}^{t}\tilde{b}_{s}^{\delta }\left( 
\tilde{\vartheta}_{s}^{\delta },X_{s}^{\delta }\right) ds+\int_{0}^{t}\sigma
_{s}\left( X_{s}^{\delta }\right) dB_{s} 
\]
where the coefficients $\tilde{b}_{s}$ are those specified in (\ref{c1}), $%
B=\left( B_{t}\right) _{t\geq 0}$ be a $d$-dimensional Brownian motion and $%
X_{0}^{\delta }\;$has distribution $\bar{\vartheta}_{0}$. Then $\tilde{%
\vartheta}_{t}^{\delta }$ will be equal to $\bar{\vartheta}_{t}^{\delta }$.
\end{proposition}

\proof%
%
Let us apply It\^{o}'s formula to the equation (\ref{new}) for $\varphi \in
C_{k}^{2}\left( {\Bbb R}^{d}\right) $. We get that 
\[
\varphi \left( X_{t}^{\delta }\right) =\varphi \left( X_{0}^{\delta }\right)
+\int_{0}^{t}\left( L_{_{s}}\varphi \left( X_{s}^{\delta }\right) +\bar{L}%
_{_{\tilde{\vartheta}_{s}^{\delta }}}^{\alpha _{s}^{\delta }}\varphi \left(
X_{s}^{\delta }\right) \right) ds+\int_{0}^{t}\sum_{i,j=1}^{d}\sigma
_{s}^{ij}\left( X_{s}\right) \partial _{x_{i}}\varphi \left( X_{s}\right)
dB_{s}^{j}, 
\]
which yields, by taking expectation and applying Proposition \ref{parts},
that $\tilde{\vartheta}^{\delta }\;$satisfies 
\begin{eqnarray*}
\tilde{\vartheta}_{t}^{\delta }\left( \varphi \right) &=&\tilde{\vartheta}%
_{0}^{\delta }\left( \varphi \right) +\int_{0}^{t}\tilde{\vartheta}%
_{s}^{\delta }\left( \left( L_{_{s}}+\bar{L}_{_{\tilde{\vartheta}%
_{s}^{\delta }}}^{\alpha _{s}^{\delta }}\right) \varphi \right) ds \\
&=&\tilde{\vartheta}_{0}^{\delta }(\varphi )+\int_{0}^{t}\tilde{\vartheta}%
_{s}^{\delta }\left( \alpha _{s}^{\delta }\varphi +L_{s}\varphi \right) -%
\tilde{\vartheta}_{s}^{\delta }\left( \alpha _{s}^{\delta }\right) \tilde{%
\vartheta}_{s}^{\delta }\left( \varphi \right) ds
\end{eqnarray*}
So $\tilde{\vartheta}^{\delta }$ satisfied the PDE (\ref{PDE2}) and, using
the uniqueness of the solution of (\ref{PDE2}), it follows that $\tilde{%
\vartheta}^{\delta }=$ $\bar{\vartheta}^{\delta }$.%
\endproof%
%

\section{ Uniqueness of the Solution}

In the following we will prove the uniqueness of a solution of (\ref{new}) \
in the class of good of processes as defined in the previous section.

\begin{theorem}
\label{unique}There exists at most one solution of (\ref{new}) which is a
good process.
\end{theorem}

\proof%
%
First, we note that if $X^{\delta }=\left( X_{t}^{\delta }\right) _{t\geq 0}$
is a solution of (\ref{new}) then the distribution of $X_{t}^{\delta }\,$
satisfies (\ref{PDE2}) and therefore is uniquely determined and equal to $%
\bar{\vartheta}_{t}^{\delta }$. Therefore we only need to justify the
uniqueness of the solution of 
\begin{equation}
dX_{t}^{i,\delta }=\tilde{b}_{t}^{\delta }\left( \bar{\vartheta}_{t}^{\delta
},X_{t}^{i,\delta }\right) dt+\sigma _{t}\left( X_{t}^{i}\right) dB_{t}^{i},
\label{new2}
\end{equation}
which is obtained from (\ref{new}) by replacing $\tilde{\vartheta}%
_{t}^{\delta }$ with $\bar{\vartheta}_{t}^{\delta }$. Hence (\ref{new2}) is
an ordinary SDE whose coefficients are locally Lipschitz. In particular $%
x\rightarrow \frac{1}{\bar{\vartheta}_{t}\left( x\right) }$ is locally
Lipschitz. Hence, for example by Theorem 3.1 page 164 in \cite{iw}, there
exists a stopping time $\zeta $ such that equation (\ref{new}) has a unique
solution in the interval $[0,\zeta )$ and, on the event, $\zeta <\infty $ we
have 
\[
\lim \sup_{t\rightarrow \zeta }\left| \left| X_{t}^{\delta }\right| \right|
=\infty \,.
\]
We want to prove next that the event $\zeta <\infty $ has null probability
(hence the solution is unique for all \ $t>0$). To do this, it is enough to
show that, for all $t>0$ we have $\left| \left| X_{t}^{\delta }\right|
\right| <\infty \,\,P-$almost surely. This fact implies, by the continuity
of the trajectories, that $\sup_{\left[ 0,t\right] }\left| \left|
X_{s}^{\delta }\right| \right| <\infty $ almost surely and hence that $\zeta
\geq t$ almost surely. For this it suffices to prove that 
\begin{equation}
\lim_{k\rightarrow \infty }P\left( \left| \left| X_{t}^{\delta }\right|
\right| \geq k\right) =0.  \label{lastcondition}
\end{equation}
We have 
\[
P\left( \left| \left| X_{t}^{\delta }\right| \right| \geq k\right) =\tilde{%
\vartheta}_{t}^{\delta }(\varphi _{k})=\bar{\vartheta}_{t}^{\delta }(\varphi
_{k})=\frac{\vartheta _{t}^{\delta }(\varphi _{k})}{\vartheta _{t}^{\delta
}(1)}\leq \frac{\vartheta _{t}^{\delta }(\bar{\varphi}_{k})}{\vartheta
_{t}^{\delta }(1)},
\]
where 
\[
\varphi _{k}\left( x\right) \triangleq I_{B\left( 0,k\right) }\left(
x\right) =\left\{ 
\begin{array}{cc}
1 & \text{if}\,\,\,\left| \left| x\right| \right| \geq k \\ 
0 & \text{otherwise.}
\end{array}
\right. ,\bar{\varphi}_{k}\left( x\right) \triangleq \left\{ 
\begin{array}{cc}
1 & \text{if}\,\,\,\left| \left| x\right| \right| \geq k \\ 
\exp \left( \frac{\left| \left| x\right| \right| -k}{\left| \left| x\right|
\right| -\frac{k}{2}}\right)  & \text{if}\,\,\,\frac{k}{2}<\left| \left|
x\right| \right| <k \\ 
0 & \left| \left| x\right| \right| <\frac{k}{2}
\end{array}
\right. 
\]
So (\ref{lastcondition}) is implied by 
\begin{equation}
\lim_{k\rightarrow \infty }\vartheta _{t}^{\delta }(\bar{\varphi}_{k})=0.
\label{lastcondition2}
\end{equation}
Now since the coefficients of the operator $L_{s}$ are all uniformly bounded
it follows that $\bar{\varphi}_{k}\in \bigcap_{s\geq 0}{\cal D}\left(
L_{s}\right) $ and that $\lim_{k\rightarrow \infty }\sup_{s\in \lbrack
0,\infty )}\left| \left| L_{s}\bar{\varphi}_{k}\right| \right| =0.\,$For
arbitrary \ $T>0$, define now the function $\Psi :\left[ 0,T\right]
\rightarrow {\Bbb R}_{+}$ 
\[
\Psi \left( t\right) =\sup_{s\in \left[ 0,t\right] }\vartheta _{t}^{\delta }(%
\bar{\varphi}_{k}).
\]
From (\ref{PDE2}) one deduces that, for arbitrary \ $T>0$ and $\ t\in \left[
0,T\right] $ we have

\[
\Psi \left( t\right) \leq \bar{\vartheta}_{0}(\bar{\varphi}_{k})+\sup_{s\in
\lbrack 0,T]}\left| \left| L_{s}\bar{\varphi}_{k}\right| \right|
\int_{0}^{T}\vartheta _{s}^{\delta }\left( 1\right) ds+K_{T,\alpha ,\beta
,\omega }^{\delta }\int_{0}^{t}\Psi \left( s\right) ds,
\]
where 
\begin{eqnarray*}
K_{T,\alpha ,\beta ,\omega }^{\delta } &=&\sup_{s\in \lbrack 0,T]}\left|
\left| \alpha _{s}\right| \right|  \\
&&+\frac{1}{2}\sum_{j=1}^{m}(\sup_{s\in \lbrack 0,T]}\left| \left| \beta
_{s}^{j}\right| \right| ^{2}+\frac{1}{\delta }\sup_{s\in \lbrack 0,T]}\left|
\left| \beta _{s}^{j}\right| \right| \max_{i=0,1,...\left[ \frac{T}{\delta }%
\right] +1}(W_{t_{i+1}}^{j}\left( \omega \right) -W_{t_{i}}^{j}\left( \omega
\right) ))
\end{eqnarray*}
Hence, by Gronwall's inequality, 
\[
\sup_{s\in \left[ 0,t\right] }\vartheta _{t}^{\delta }(\bar{\varphi}%
_{k})\leq e^{K_{T,\alpha ,\beta ,\omega }^{\delta }T}\left( \bar{\vartheta}%
_{0}(\bar{\varphi}_{k})+\sup_{s\in \lbrack 0,T]}\left| \left| L_{s}\bar{%
\varphi}_{k}\right| \right| \int_{0}^{T}\vartheta _{s}^{\delta }\left(
1\right) ds\right) 
\]
which implies that 
\begin{equation}
\lim_{k\rightarrow \infty }\sup_{s\in \left[ 0,t\right] }\vartheta
_{t}^{\delta }(\bar{\varphi}_{k})=0.  \label{lastcondition3}
\end{equation}
as $\lim_{k\rightarrow \infty }\bar{\vartheta}_{0}(\bar{\varphi}_{k})=0$.
But (\ref{lastcondition3}) implies (\ref{lastcondition2}), which in turn
implies (\ref{lastcondition}) and that completes the proof of the Theorem.%
\endproof%
%

\section{Existence of the Solution}

We are now ready to complete the last step of the programme. For this we
need to add one final condition:

{\bf ES}: Let $\left( s,x\right) \longrightarrow \tilde{b}_{s}^{\delta ,\bar{%
\vartheta}_{s}^{\delta }}\left( x\right) \;$be the function defined in (\ref
{c1}) where the measure $\mu \;$is chosen to be $\bar{\vartheta}_{s}^{\delta
}$.$\;$In other words, 
\begin{equation}
\tilde{b}_{s}^{\delta ,\bar{\vartheta}_{s}^{\delta }}\left( x\right)
\triangleq b_{s}\left( x\right) +\frac{\Lambda _{\bar{\vartheta}_{s}^{\delta
}}\alpha _{s}^{\delta }\left( x\right) }{\bar{\vartheta}_{s}^{\delta }\left(
x\right) }  \label{c2}
\end{equation}
We will assume that $\left( s,x\right) \longrightarrow \tilde{b}_{s}^{\delta
,\bar{\vartheta}_{s}^{\delta }}\left( x\right) \;$is globally Lipschitz.\
More precisely we will assume that, for any $T\geq 0$. there exists a
constant $K_{T}\;$such that 
\begin{equation}
\left| \left| \tilde{b}_{s}^{\delta ,\bar{\vartheta}_{s}^{\delta }}\left(
x\right) -\tilde{b}_{s}^{\delta ,\bar{\vartheta}_{s}^{\delta }}\left(
y\right) \right| \right| \leq K_{T}\left| \left| x-y\right| \right|
\label{Lip1}
\end{equation}
for all $x,y\in R^{d}\;$and\ $s\in \left[ 0,T\right] .$

\begin{theorem}
\label{exist}Under the conditions {\bf EU+PD+ES}, equation (\ref{new}) has a
solution which is a good process.
\end{theorem}

\proof%
%
We prove the existence of the solution on an arbitrary time interval $\left[
0,T\right] .\;$From\ (\ref{Lip1}) and the fact that the function 
\[
s\in \left[ 0,T\right] \longrightarrow b_{s}\left( 0\right) +\frac{\Lambda _{%
\bar{\vartheta}_{s}^{\delta }}\alpha _{s}^{\delta }\left( 0\right) }{\bar{%
\vartheta}_{s}^{\delta }\left( 0\right) }
\]
is continuous, hence bounded, it follows (for example, by using Theorem 2.9,
p. 289, in \cite{ks}) that the equation 
\begin{equation}
d\tilde{X}_{t}^{i,\delta }=\tilde{b}_{s}^{\delta ,\bar{\vartheta}%
_{s}^{\delta }}\left( \tilde{X}_{t}^{i,\delta }\right) dt+\sigma _{t}\left( 
\tilde{X}_{t}^{i}\right) dB_{t}^{i},  \label{new3}
\end{equation}
has a (unique) solution whose distribution satisfies the linear PDE 
\begin{equation}
\tilde{\vartheta}_{t}^{\delta }(\varphi )=\bar{\vartheta}_{0}(\varphi
)+\int_{0}^{t}\tilde{\vartheta}_{s}^{\delta }\left( L_{_{s}}\varphi +\bar{L}%
_{_{\bar{\vartheta}_{s}^{\delta }}}^{\alpha _{s}^{\delta }}\varphi \right) ds
\label{PDE3}
\end{equation}
which has a unique solution (see ... in \cite{lsu}). From (\ref{PDE1}) it
follows that $\bar{\vartheta}_{s}^{\delta }\;$is a solution of PDE (\ref
{PDE3}), hence\ $\tilde{\vartheta}^{\delta }$ and $\bar{\vartheta}^{\delta },
$hence $\tilde{X}^{i,\delta }\;$is in fact a solution of the nonlinear SDE (%
\ref{new}).%
\endproof%
%

\begin{remark}
If $d=1$\ and for all $s\geq 0,\;$the\ density of$\;\vartheta _{s}^{\delta }$
with respect to the Lebesgue\ measure has the form 
\[
\vartheta _{s}^{\delta }\left( x\right) =e^{-F_{s}^{\delta }\left( x\right)
},\;x\in {\Bbb R}. 
\]
where $F_{s}^{\delta }\;$is a differentiable convex function, condition {\bf %
ES\ }is satisfied.
\end{remark}

\proof%
%
We have three cases:\newline
{\bf 1.} If $\frac{dF_{s}^{\delta }}{dx}\left( x\right) =0,$ then $\frac{d}{%
dx}\left( \frac{\Lambda _{\vartheta _{s}^{\delta }}\alpha _{s}^{\delta
}\left( x\right) }{e^{-F_{s}^{\delta }\left( x\right) }}\right) =\left(
\alpha _{s}^{\delta }\right) ^{\mu }\left( x\right) \;$which is bounded
since $\alpha _{s}^{\delta }$ is bounded.\newline
{\bf 2. }If $\frac{dF_{s}^{\delta }}{dx}\left( x\right) >0,$ then using the
first part of the representation (\ref{onerep}) we get 
\[
\frac{d}{dx}\left( \frac{\Lambda _{\vartheta _{s}^{\delta }}f\left( x\right) 
}{e^{-F_{s}^{\delta }\left( x\right) }}\right) =\left( \alpha _{s}^{\delta
}\right) ^{\mu }\left( x\right) +\int_{x}^{\infty }\left( \alpha
_{s}^{\delta }\right) ^{\mu }\left( y\right) \frac{dF_{s}^{\delta }}{dx}%
\left( x\right) e^{F_{s}^{\delta }\left( x\right) -F_{s}^{\delta }\left(
y\right) }dy
\]
Since $F_{s}^{\delta }$\thinspace is convex, for $y>x\;$we have 
\[
\frac{F_{s}^{\delta }\left( y\right) -F_{s}^{\delta }\left( x\right) }{y-x}%
\geq \frac{dF_{s}^{\delta }}{dx}\left( x\right) \Longrightarrow
F_{s}^{\delta }\left( x\right) -F_{s}^{\delta }\left( y\right) \leq \left(
x-y\right) \frac{dF_{s}^{\delta }}{dx}\left( x\right) 
\]
So 
\begin{eqnarray*}
\left| \int_{x}^{\infty }\left( \alpha _{s}^{\delta }\right) ^{\mu }\left(
y\right) \frac{dF_{s}^{\delta }}{dx}\left( x\right) e^{F_{s}^{\delta }\left(
x\right) -F_{s}^{\delta }\left( y\right) }dy\right|  &\leq &\left| \left|
\left( \alpha _{s}^{\delta }\right) ^{\mu }\right| \right| \frac{%
dF_{s}^{\delta }}{dx}\left( x\right) \int_{x}^{\infty }e^{\left( x-y\right) 
\frac{dF_{s}^{\delta }}{dx}\left( x\right) }dy \\
&\leq &\left| \left| \left( \alpha _{s}^{\delta }\right) ^{\mu }\right|
\right| 
\end{eqnarray*}
since $\frac{dF_{s}^{\delta }}{dx}\left( x\right) \int_{x}^{\infty
}e^{\left( x-y\right) \frac{dF_{s}^{\delta }}{dx}\left( x\right) }dy=1.$%
\newline
{\bf 3. }If{\bf \ }$\frac{dF_{s}^{\delta }}{dx}\left( x\right) <0,$then then
using the second part of the representation (\ref{onerep}) we get 
\[
\frac{d}{dx}\left( \frac{\Lambda _{\vartheta _{s}^{\delta }}f\left( x\right) 
}{e^{-F_{s}^{\delta }\left( x\right) }}\right) =\left( \alpha _{s}^{\delta
}\right) ^{\mu }\left( x\right) -\int_{-\infty }^{x}\left( \alpha
_{s}^{\delta }\right) ^{\mu }\left( y\right) \frac{dF_{s}^{\delta }}{dx}%
\left( x\right) e^{F_{s}^{\delta }\left( x\right) -F_{s}^{\delta }\left(
y\right) }dy
\]
and we follow the same steps as in the previous case.%
\endproof%
%

We complete the section by noting that the above construction works with
minimal changes when (the initial Markov process) $\xi \;$is a reflecting
boundary diffusion. In this case, the analysis simplifies considerably if
the domain is compact. For example, the cumbersome condition {\bf ES} is
replaced by the assumption that $\xi \,$has a density which is bounded away
from 0. We will detail the analysis of this case in a forthcoming paper,
together with the description of the associated numerical algorithm.
However, for completeness, we briefly present the results here:

Following the notation and results in \cite{ls}, assume that $\xi \;$is the
solution of {\em a stochastic differential equation with reflection along
the normal}. In other words, $\xi \;$is the unique solution of the equation 
\begin{equation}
\xi _{t}=\xi _{0}+\int_{0}^{t}b_{s}\left( \xi _{s}\right)
ds+\int_{0}^{t}\sigma _{s}\left( \xi _{s}\right) dB_{s}-k_{t},  \label{sder}
\end{equation}
where $\xi _{t}\in $ $\overline{B\left( 0,M\right) }$\footnote{$B\left(
0,M\right) $\ is the ball of center $0$ and radius $M.$} for all $t\geq 0$
and $k\;$is a bounded variation process 
\[
\left| k\right| _{t}=\int_{0}^{t}1_{\left\{ \xi _{s}\in \partial B\left(
0,M\right) \right\} }d\left| k\right| _{s},\,\,\,k_{t}=\int_{0}^{t}\frac{\xi
_{s}}{M}d\left| k\right| _{s}
\]
For general domains $D,\;k$\ is defined as 
\[
\left| k\right| _{t}=\int_{0}^{t}1_{\left\{ \xi _{s}\in \partial D\right\}
}d\left| k\right| _{s},\,\,\,k_{t}=\int_{0}^{t}n\left( \xi _{s}\right)
d\left| k\right| _{s}
\]
where $n\left( x\right) $ is the unit outward normal to $\partial D$ at $x$. 

Obviously the generator  $L=\left( L_{s}\right) _{s\geq 0}$ associated to $%
\xi $ has the form (\ref{operator}) for any $\varphi $ in $C_{K}^{2}\left(
B\left( 0,M\right) \right) ,\;$the set of twice differentiable functions
defined on $B\left( 0,M\right) $ with compact support. The operator $\Lambda
_{\mu }f$ \ is\ well defined for measures $\mu \;$with support in $\overline{%
B\left( 0,M\right) }\;$absolute continuous with respect to the Lebesgue
measure and with bounded density. The definition of a {\it good} process is
now slightly simpler: A continuous process $X=\left( X_{t}\right) _{t\geq
0}\;$which takes values in $\overline{B\left( 0,M\right) }$ is a {\em good}
process if the following two conditions are satisfied ( as before, $\tilde{%
\vartheta}_{t}$ is the distribution of $X_{t}$):\bigskip \newline
$\blacklozenge $\ $\tilde{\vartheta}_{t}$ is absolutely continuous with
respect to the Lebesgue measure and its density $x\rightarrow \tilde{%
\vartheta}_{t}\left( x\right) $ is positive on $\overline{B\left( 0,M\right) 
}$.\newline
$\blacklozenge $\ The function \ $\left( t,x\right) \in \left[ 0,\infty %
\right] \times $ $\overline{B\left( 0,M\right) }\longrightarrow \tilde{%
\vartheta}_{t}\left( x\right) \,\,$is continuous and bounded.\medskip 
\newline
Proposition (\ref{nonproc}) now becomes:

\begin{proposition}
Let $\left( \tilde{\Omega},{\cal \tilde{F}},\tilde{P}\right) \;$be a
probability space on which there exists a good process $X^{\delta }=\left(
X_{t}^{\delta }\right) _{t\geq 0}$ which satisfies equation 
\begin{equation}
X_{t}^{\delta }=X_{0}^{\delta }+\int_{0}^{t}\tilde{b}_{s}^{\delta }\left( 
\tilde{\vartheta}_{s}^{\delta },X_{s}^{\delta }\right) ds+\int_{0}^{t}\sigma
_{s}\left( X_{s}^{\delta }\right) dB_{s}-k_{t}^{\delta }  \label{new4}
\end{equation}
where $X_{t}^{\delta }\in $ $\overline{B\left( 0,M\right) }$ for all $t\geq
0 $ and $k^{\delta }\;$is a bounded variation process 
\[
\left| k^{\delta }\right| _{t}=\int_{0}^{t}1_{\left\{ X_{s}^{\delta }\in
\partial B\left( 0,M\right) \right\} }d\left| k^{\delta }\right|
_{s},\,\,k_{t}^{\delta }=\int_{0}^{t}\frac{X_{s}^{\delta }}{M}d\left|
k^{\delta }\right| _{s} 
\]
In (\ref{new4}),\ the coefficients $\tilde{b}_{s}$ are those specified in (%
\ref{c1}), $B=\left( B_{t}\right) _{t\geq 0}$ be a $d$-dimensional Brownian
motion and $X_{0}^{\delta }\;$has distribution $\bar{\vartheta}_{0}$. Then
the distribution of $X_{t}^{\delta }$ will be equal to $\bar{\vartheta}%
_{t}^{\delta }$.
\end{proposition}

The uniqueness of the pair $\left( X^{\delta },\,k^{\delta }\right) \;$%
follows by slight variation of the argument for the proof of theorem (\ref
{unique}). The existence of the pair $\left( X^{\delta },\,k^{\delta
}\right) $\ requires the following simpler condition:\medskip

{\bf ES'}:\ Assume that, $\xi =\left( \xi _{t}\right) _{t\geq 0},$\ the
solution of the SDE (\ref{sder}) is a good process and that its density $%
\breve{\vartheta}_{t}\left( x\right) \;$is\ bounded away from 0. More
precisely, for any $T>0\;$we have 
\[
\inf_{t\in \left[ 0,T\right] x\in \overline{B\left( 0,M\right) }}\breve{%
\vartheta}_{t}\left( x\right) >0. 
\]
Then we have the following result:

\begin{theorem}
Under the conditions {\bf EU+PD+ES'}, equation (\ref{new4}) has a unique
solution which is a good process.
\end{theorem}

The results follows by proving that {\bf ES'\ }implies{\bf \ ES}. Since $%
\bar{\vartheta}_{t}^{\delta }\;$has the representation 
\begin{equation}
\bar{\vartheta}_{t}^{\delta }\left( \varphi \right) =E\left[ \varphi \left(
\xi _{t}\right) \exp \left( \int_{0}^{t}\alpha _{s}^{\delta }\left( \xi
_{s}\right) -\bar{\vartheta}_{t}^{\delta }\left( \alpha _{s}^{\delta
}\right) ds\right) \right] .  \label{fkp2}
\end{equation}
$\bar{\vartheta}_{t}^{\delta }\;$is\ absolutely continuous with respect to
the Lebesgue measure, too. Further, since $\alpha _{s}^{\delta }\;$as
defined in (\ref{alpha}) is a bounded function on $\overline{B\left(
0,M\right) },$\ we deduce that there exists a positive constant $\epsilon \;$%
such that 
\[
\bar{\vartheta}_{t}^{\delta }\left( \varphi \right) \geq \epsilon E\left[
\varphi \left( \xi _{t}\right) \right] =\epsilon \int_{B\left( 0,M\right)
}\varphi \left( x\right) \breve{\vartheta}_{t}\left( x\right) dx
\]
Hence the density of $\bar{\vartheta}_{t}^{\delta }\;$is uniformly bounded
away from 0. One can also prove that the density of $\bar{\vartheta}%
_{t}^{\delta }\;$has uniformly bounded first order derivatives. This implies
that the function $\left( s,x\right) \longrightarrow \frac{\Lambda _{\bar{%
\vartheta}_{s}^{\delta }}\alpha _{s}^{\delta }\left( x\right) }{\bar{%
\vartheta}_{s}^{\delta }\left( x\right) }\;$is globally Lipschitz, hence 
{\bf ES\ }is satisfied.

\begin{remark}
Numerical methods for equations of type (\ref{sder}) have extensively
developed (see for example \cite{pe}, \cite{slo} and the references
therein).\ We will adapt these schemes in order to produce a numerical
method for solving (\ref{new4}).
\end{remark}

\section{Rates of convergence and final remarks}

It is obvious that, for any $\varphi \in C_{b}\left( {\Bbb R}^{d}\right) ,\;$%
we have \ 
\[
\tilde{E}\left[ \left( \bar{\vartheta}_{t}^{\delta }\left( \varphi \right) -%
\bar{\vartheta}_{t}^{\delta ,n}\left( \varphi \right) \right) ^{4}\right]
\leq \frac{\left| \left| \varphi \right| \right| ^{4}}{n^{2}}
\]
so 
\[
\tilde{E}\left[ \left| \left| \bar{\vartheta}_{t}^{\delta }-\bar{\vartheta}%
_{t}^{\delta ,n}\right| \right| _{{\cal M}}^{4}\right] \leq \frac{1}{n^{2}}.
\]
where the norm $\left| \left| \cdot \right| \right| _{{\cal M}}$\ is defined
in (\ref{mnorm}). From (\ref{rate2}), we have that 
\[
\left| \left| \vartheta _{t}^{\delta \left( \omega ,n\right) }-\vartheta
_{t}\right| \right| _{{\cal M}}\leq \frac{K}{\sqrt{n}}\text{ and}\;\;\left|
\left| \bar{\vartheta}_{t}^{\delta \left( \omega ,n\right) }-\bar{\vartheta}%
_{t}\right| \right| _{{\cal M}}\leq \frac{\bar{K}}{\sqrt{n}},
\]
hence 
\[
\tilde{E}\left[ \left| \left| \bar{\vartheta}_{t}-\bar{\vartheta}%
_{t}^{\delta \left( \omega ,n\right) ,n}\right| \right| _{{\cal M}}^{4}%
\right] \leq \frac{\left( \sqrt[4]{\bar{K}}+1\right) ^{4}}{n^{2}},
\]
so $\bar{\vartheta}_{t}^{\delta \left( \omega ,n\right) ,n}\;$converges\ to $%
\bar{\vartheta}_{t}\;\;\tilde{P}$-almost surely. We also have almost sure
convergence if we view $\bar{\vartheta}_{t}^{\delta \left( \omega ,n\right)
,n}$\ and $\bar{\vartheta}_{t}\;$as processes on the product space $(\bar{%
\Omega},\bar{{\cal F}},\bar{P})$ 
\[
(\bar{\Omega},\bar{{\cal F}},\bar{P})=(\Omega \times \hat{\Omega},{\cal F}%
\tbigotimes {\cal F},\tilde{P}\tbigotimes \tilde{P})
\]
on which we `lift' the processes $\bar{\vartheta}_{t}^{\delta \left( \omega
,n\right) ,n}$ and $\bar{\vartheta}_{t}$ from the component spaces. Finally,
in introduction we chose the asymptotic representation for $\vartheta _{t}\;$%
to be
\begin{equation}
\vartheta _{t}=\lim_{n\rightarrow \infty }\frac{a_{t}^{n}}{n}%
\sum_{i=1}^{n}\delta _{X_{t}^{i,\delta (W_{\cdot }\left( \omega \right)
,n)}},  \label{aoleo}
\end{equation}
where 
\begin{eqnarray}
a_{t}^{n} &=&\exp \left( \int_{0}^{t}\bar{\vartheta}_{t}^{\delta \left(
\omega ,n\right) ,n}\left( \alpha _{s}\right) ds+\sum_{j=1}^{m}\int_{0}^{t}%
\bar{\vartheta}_{t}^{\delta \left( \omega ,n\right) ,n}\left( \beta
_{s}^{j}\right) d\left( W_{s}^{\delta \left( \omega ,n\right) }\right)
^{j}\right.   \nonumber \\
&&\left. \hspace{1in}-\frac{1}{2}\sum_{j=1}^{m}\int_{0}^{t}\bar{\vartheta}%
_{t}^{\delta \left( \omega ,n\right) ,n}\left( \beta _{s}^{j}\right)
^{2}ds\right)   \label{weightan}
\end{eqnarray}
Then each of the terms in the formula (\ref{weightan}), i.e., $\int_{0}^{t}%
\bar{\vartheta}_{t}^{\delta \left( \omega ,n\right) ,n}\left( \alpha
_{s}\right) ds,$\newline
$\int_{0}^{t}\bar{\vartheta}_{t}^{\delta \left( \omega ,n\right) ,n}\left(
\beta _{s}^{j}\right) d\left( W_{s}^{\delta \left( \omega ,n\right) }\right)
^{j}$ and $\int_{0}^{t}\bar{\vartheta}_{t}^{\delta \left( \omega ,n\right)
,n}\left( \beta _{s}^{j}\right) ^{2}ds\;$converge $\bar{P}-$almost surely to 
$\int_{0}^{t}\bar{\vartheta}_{t}\left( \alpha _{s}\right) ds,\;\int_{0}^{t}%
\bar{\vartheta}_{t}\left( \beta _{s}^{j}\right) dW_{s}^{j}$ and,
respectively, $\int_{0}^{t}\bar{\vartheta}_{t}\left( \beta _{s}^{j}\right)
^{2}ds,$ hence $a_{t}^{n}\;$converges almost surely to $\vartheta _{t}\left(
1\right) $, so (\ref{weightan}) holds true.


\begin{thebibliography}{99}
\bibitem{c}  {\sc J. M. C. Clark, }{\em The design of robust approximations
to the stochastic differential equations of nonlinear filtering}.
Communication systems and random process theory (Proc. 2nd NATO Advanced
Study Inst., Darlington, 1977), pp. 721--734. NATO Advanced Study Inst.
Ser., Ser. E: Appl. Sci., No. 25, Sijthoff \& Noordhoff, Alphen aan den
Rijn, 1978.

\bibitem{cc}  {\sc J. M. C. Clark, D. Crisan}, {\em On a robust version of
the integral representation formula of nonlinear filtering}, to appear in
Probability Theory and Related Fields, 2005.

\bibitem{dan}  {\sc D. Crisan}, {\em Numerical methods for solving the
stochastic filtering problem}, Numerical methods and stochastics (Toronto,
ON, 1999), 1--20, Fields Inst. Commun., 34, Amer. Math. Soc., Providence,
RI, 2002.

\bibitem{dan2}  {\sc D. Crisan, J. Gaines, T. Lyons, } {\em Convergence of a
branching particle method to the solution of the Zakai equation}. SIAM J.
Appl. Math. 58 (1998), no. 5, 1568--1590 (electronic).

\bibitem{dan3}  {\sc D. Crisan, T. Lyons,} {\em A particle approximation of
the solution of the Kushner-Stratonovitch equation} Probab. Theory Related
Fields 115 (1999), no. 4, 549--578.

\bibitem{dan4}  {\sc D. Crisan, P. Del Moral T. Lyons, }{\em Interacting
particle systems approximations of the Kushner-Stratonovitch equation.} Adv.
in Appl. Probab. 31 (1999), no. 3, 819--838.

\bibitem{dan5}  {\sc D. Crisan,} Particle filters---a theoretical
perspective. Sequential Monte Carlo methods in practice, 17--41, Stat. Eng.
Inf. Sci., Springer, New York, 2001.

\bibitem{dan6}  {\sc D. Crisan, A. Doucet}, {\em A survey of convergence
results on particle filtering methods for practitioners.} IEEE Trans. Signal
Process. 50, no. 3, pp 736--746, 2002.

\bibitem{dan7}  {\sc D. Crisan,} {\em Exact rates of convergence for a
branching particle approximation to the solution of the Zakai equation }Ann.
Probab. 31 (2003), no. 2, 693--718.

\bibitem{dmm}  {\sc P. Del Moral, L. Miclo,} {\em Branching and interacting
particle systems approximations of Feynman-Kac formulae with applications to
non-linear filtering.} S\'{e}minaire de Probabilit\'{e}s, XXXIV, 1--145,
Lecture Notes in Math., 1729, Springer, Berlin, 2000.

\bibitem{d2}  {\sc T. E. Duncan}, {\em Likelihood functions for stochastic
signals in white noise}, Information and Control 16 1970 303--310.

\bibitem{fkk}  {\sc M. Fujisaki, G. Kallianpur, H. Kunita,} {\em Stochastic
differential equations for the non linear filtering problem}. Osaka J. Math.
9, pp 19--40, 1972.

\bibitem{iw}  {\sc N. Ikeda, S. Watanabe,} {\em Stochastic differential
equations and diffusion processes,} North-Holland Mathematical Library, 24.
North-Holland Publishing Co., Amsterdam-New York; Kodansha, Ltd., Tokyo,
1981.

\bibitem{ks}  {\sc I. Karatzas, S. E. Shreve,} {\em Brownian motion and
stochastic calculus,} Second Edition. Graduate Texts in Mathematics, 113.
Springer-Verlag, New York, 1991.

\bibitem{kj1}  {\sc T. G. Kurtz, J. Xiong}, {\em Particle Representation for
a Class of Nonlinear SPDE's}, Stochastic Process. Appl. 83, no. 1, pp
103--126, 1999.

\bibitem{kj2}  {\sc T. G. Kurtz, J. Xiong}, {\em Numerical Solutions for a
Class of SPDEs with Application to Filtering.} Stochastics in finite and
infinite dimensions, 233--258, Trends Math., Birkh\"{a}user Boston, 2001.

\bibitem{ku}  {\sc H. J. Kushner, }{\em Dynamical equations for optimal
nonlinear filtering}, J. Differential Equations 3, 179--190, 1967.

\bibitem{lsu}  {\sc O. A. Lady\~{z}enskaja, V. A. Solonnikov, N. N.
Ural'ceva, } {\em Linear and quasilinear equations of parabolic type}.
(Russian) Translated from the Russian by S. Smith. Translations of
Mathematical Monographs, Vol. 23 American Mathematical Society, Providence,
R.I. 1967.

\bibitem{ls}  {\sc P. L. Lions, A.S. Sznitman,} {\em Stochastic differential
equations with reflecting boundary conditions}, Comm. Pure Appl. Math. 37,
no. 4, 511--537, 1984.

\bibitem{liptser}  {\sc R. S. Liptser, A. N. Shiryaev}, {\em Statistics of
random processes. I. General theory,} Second, revised and expanded edition,
Springer-Verlag, Berlin, 2001.

\bibitem{m1}  {\sc R. E. Mortensen,} {\em Stochastic optimal control with
noisy observations}, Internat. J. Control (1) 4 1966 455--464.

\bibitem{Pardoux}  {\sc E. Pardoux}, {\em Filtrage Non Lin\'{e}aire et
Equations aux D\'{e}riv\'{e}es Partielles Stochastiques Associ\'{e}es},
Ecole d'Et\'{e} de Probabilit\'{e}s de Saint-Flour XIX - 1989, Lecture Notes
in Mathematics, 1464, Springer-Verlag, 1991.

\bibitem{pe}  {\sc R. Pettersson}, {\em Penalization schemes for reflecting
stochastic differential equations}, Bernoulli 3, no. 4, 403--414, 1997.

\bibitem{slo}  {\sc L. Slominski}, {\em Euler's approximations of solutions
of SDEs with reflecting boundary}, Stochastic Process. Appl. 94, no. 2,
317--337, 2001.

\bibitem{roz}  {\sc B. L. Rozovskii}, {\em Stochastic Evolution Systems, }%
Kluwer, Dordrecht, 1990.

\bibitem{lsc}  {\sc L. Saloff-Coste,} {\em Aspects of Sobolev-type
inequalities,} London Mathematical Society Lecture Note Series, 289.
Cambridge University Press, Cambridge, 2002.

\bibitem{za}  {\sc M Zakai, }{\em On the optimal filtering of diffusion
processes}, Z. Wahrscheinlichkeitstheorie und Verw. Gebiete, 11, pp
230--243, 1969.
\end{thebibliography}
\end{document}